\documentclass[a4paper, 11pt]{article}
\usepackage{mathrsfs}
\usepackage{amsfonts}
\usepackage{amsmath,amsthm,mathrsfs}
\usepackage{amsfonts,amssymb}

\begin{document}
\title{\bf Wavelet transform and Radon transform\\ on the quaternion Heisenberg group}

\author{\small Jianxun He \footnote{Supported by the National Natural Science Foundation of China \#10971039
and the Specialized Research Fund for the  Doctoral Program of Higher Education of China
under Grant \#200810780002. }\\[-0.1cm]
\small School of Mathematics and Information Sciences\\[-0.1cm]
\small Guangzhou University\\[-0.1cm]
\small Guangzhou 510006, P.R. China\\[-0.1cm]
\small E-mail: hejianxun@gzhu.edu.cn  \\[0.15cm]
\small Heping Liu \footnote{Supported by the National Natural Science Foundation of
China under Grant \#10871003, \#10990012 and the Specialized
Research Fund for the Doctoral Program of Higher Education of China
under Grant \#2007001040.}\\[-0.1cm]
\small LMAM, School of Mathematical Sciences\\[-0.1cm]
\small Peking University\\[-0.1cm]
\small Beijing 100871, P.R. China\\[-0.1cm]
\small E-mail: hpliu@pku.edu.cn \\[0.15cm]}

\date{}
\maketitle

\begin{abstract}

Let $\mathscr Q$ be the quaternion Heisenberg group, and let
$\mathbf P$ be the affine automorphism group of $\mathscr Q$. We
develop the theory of continuous wavelet transform on the quaternion
Heisenberg group via the unitary representations of $\mathbf P$ on
$L^2(\mathscr Q)$. A class of radial wavelets is constructed. The
inverse wavelet transform is simplified by using radial wavelets.
Then we investigate the Radon transform on $\mathscr Q$. A
Semyanistri-Lizorkin
 space is introduced, on which the Radon
transform is a bijection.
We deal with the Radon transform on $\mathscr Q$ both by the Euclidean Fourier
transform and the group Fourier transform. These two treatments are essentially
equivalent. We also give an inversion formula by using wavelets, which does not
require the smoothness of functions if the wavelet is smooth. \\


\end{abstract}

\maketitle

\newcommand\sfrac[2]{{#1/#2}}
\newcommand\cont{\operatorname{cont}}
\newcommand\diff{\operatorname{diff}}

\vskip 0.5 true cm
\section{\bf Introduction}
\vskip 0.5 true cm

The Heisenberg group, denoted by $\mathbf H_n$, is the simplest
example of two-step nilpotent Lie group. There are many works
devoted to the theory of harmonic analysis on this group. Geller
\cite{Ge} established the theory of Fourier analysis on $\mathbf
H_n$. More results can be found in \cite{Fo}, \cite{Mu}, \cite{St},
\cite{Str} and the references therein. Wavelet analysis on the
Euclidean space $\mathbb R^n$ has many applications in pure and
applied mathematics (see \cite{Da}). It is important to extend the
theory of wavelet analysis to various cases. Several authors
developed the theory of continuous wavelet transform on the
Heisenberg group $\mathbf H_n$ (see \cite{JL4}, \cite{Li}).
Recently, some further extensions of wavelet analysis were published
in \cite{Ge1}, \cite{I}. The Radon transform represents an
interesting object from the point of view of both harmonic analysis
and integral geometry. Also it is a very useful tool to deal with
the problems of mathematics and engineer. The Radon transform on the
Heisenberg group (or Heisenberg Radon transform) was studied by
Geller-Stein \cite{GS} and Strichartz \cite{St}. Related further
extension we refer the reader to see \cite{Fel}, \cite{Pe}.
Holschneider (\cite{Ho}) is the first author who applied the inverse
wavelet transform to the inverse Radon transform on the
two-dimensional plane.  Rubin \cite{Ru} extended to the case of
$k$-dimensional Radon transform on $\mathbb R^n$. Nessibi-Trim\`eche
\cite{Ne} obtained an inversion formula of the Radon transform by
using generalized wavelets on the Laguerre hypergroup. The further
development can be found in \cite{Fel}, \cite{J1} $-$ \cite{JL3} and
\cite{Ru1}. Heisenberg type groups are generalizations of the
Heisenberg group, which are important both on geometry and analysis
(see \cite{Co}, \cite{Ka}, \cite{Kar} and \cite{Ko}). The quaternion
Heisenberg group is a typical Heisenberg type group other than the
Heisenberg group, which has a good explain in geometry (see
\cite{Be}, \cite{Ch}).  Tie-Wong \cite{Tie} studied the heat kernel
and Green functions associated with  the Sub-Laplacians. Moreover,
Zhu \cite{Zh} investigated the property of the Riesz transforms on
this group. In this article we study the wavelet transform and the
Radon transform on the quaternion Heisenberg group $\mathscr Q$.

This article is organized as follows. In the remainder of this
section we shall recall some basic facts of quaternion
numbers, and state the background of the quaternion Heisneberg group
$\mathscr Q$. Also, we describe the automorphism group of $\mathscr
Q$ and its regular representation. In Section 2 we give the direct sum decomposition for $L^2(\mathscr Q)$
in terms of the group Fourier transform, in which every subspace is irreducible for the representation of
the automorphism group of $\mathscr Q$. The theory of continuous wavelet transform
via square integral group representation (see \cite{Fe}) is developed in Section 3.
we also construct a class of radial wavelets. The inverse wavelet transform can be simplified by using radial wavelets.
Section 4 is devoted to the Radon transform. We introduce a Semyanistri-Lizorkin type space, on which the
Radon transform is a bijection. Then we give the inverse Radon transform. We give the results in two different ways.
One way is in terms of the Euclidean Fourier transform and another the group Fourier transform.
We also prove that these two treatments are essentially equivalent.
Finally, in Section 5 we make use of inverse wavelet transform to derive an
inversion formula of the Radon transform in $L^2$-sense,
which does not require the smoothness of functions if the wavelet is smooth.

\vskip 0.3 true cm

 Let $\mathbb Q$ denote the set of all quaternion numbers, $i,j,k$
are the three imaginary units satisfying: $i^2=j^2=k^2=ijk=-1$.
For any $q\in \mathbb Q$, we can write $q=q_0+q_1i+q_2j+q_3k$. For
convenience, we also set $q=(q_0,q_1,q_2,q_3)$ where
$q_0,q_1,q_2,q_3\in \mathbb R$. Let $\Re q$ and $\Im q$ denote the
real part and imaginary part of $q$ respectively. Then $\Re q=q_0$,
$\Im q=q_1i+q_2j+q_3k=(q_1,q_2,q_3)=q^{I}$. In contrast to
complex numbers, $\Im q$ is not a real number. The multiplication of two quaternion numbers $q, h$ is given
by
$$
\Re (qh)=q_0h_0-q^{I}\cdot h^{I}, \,\, \Im
(qh)=q_0h^{I}+h_0q^{I}+q^{I}\times h^{I}.
$$
$\overline q=q_0-q_1i-q_2j-q_3k$ denotes the conjugate of $q$. The
scalar product is given by $\langle q,h \rangle= \Re
(\overline{q}h)$. The norm is $|q|^2= \langle q, q \rangle =
\sum_{l=0}^3q_l^2$. Then we have $\overline {qh}=\overline h
\overline q,
 |qh|=|q||h|, q^{-1}=\frac {\overline q}{|q|^2}$. We also note
 the following facts. The reals are only quaternions which commute
 with all quaternions, and $q^2=-1$ if and only if $|q|=1$ and
 $q=\Im q$. We will identify $\mathbb Q$ with $\mathbb R^4$ and $\Im \mathbb Q$ with $\mathbb R^3$ if necessary.

Similar to the Heisenberg group, the quaternion Heisenberg group is the boundary of the Siegel
 upper-half space in the quaternion content ( see \cite{St}). Let $ \mathcal B$ be the unit
 ball in $\mathbb Q^2$ which is given by
 $$
  \mathcal B=\left\{(h_1,h_2):|h_1|^2+|h_2|^2<1\right\}.
 $$
Then $ \mathcal B$ is biholomorphic to the  Siegel upper-half plane
in $\mathbb Q^2$ given by
 $$
 \mathcal
 U=\left\{(q_1,q_2)\in \mathbb Q^2: \Re q_2>|q_1|^2\right\}.
 $$
 The Cayley transform is given by
$$
\left\{
\begin{array}{ll}
 q_1&=\displaystyle{\frac {h_1}{1+h_2}=\frac {h_1(1+\overline h_2)}{|1+h_2|^2}}\\[4mm]
 q_2&=\displaystyle{\frac
{1-h_2}{1+h_2}=\frac {(1-h_2)(1+\overline h_2)}{|1+h_2|^2}}.
\end{array}
\right.
$$
 Let $r=r(q_1,q_2)=\Re q_2-|q_1|^2$ be the height function.
 Setting
 $$
 x=q_1,  t=\Im q_2, r=r(q_1,q_2)=\Re q_2-|q_1|^2.
 $$
 If we adopt the Heisenberg coordinate $(x,t,r)$, then the Siegel upper-half plane is denoted by
 $$
 \mathcal U=\left\{(x,t,r):x\in \mathbb Q, t \in \Im \mathbb Q, r>0 \right\}.
 $$
 The boundary of $\mathcal U$ can be identified with the quaternion Heisenberg group
 denoted by
 $$
 \mathscr Q=\left\{(x,t): x\in \mathbb Q, t\in \Im \mathbb Q\right\}.
 $$
 The group multiplication  is given by
 $$
 (x,t)(x^{\prime},t^{\prime})=(x+x^{\prime},
 t+t^{\prime}-2\Im(\overline {x^{\prime}}x)).
 $$
 Thus we see that $\mathscr Q\cong \mathbb R^4\times \mathbb
 R^3$. The Haar measure on $\mathscr Q$ coincides  with the
 Lebesgue measure on $\mathbb R^4\times \mathbb R^3$  which is
 denoted by $dxdt$.

Let $u\in \mathbb Q$, $u=a+bi+cj+dk=(a+bi)+(c+di)j$. This
implies the identification of $\mathbb Q$ with $\mathbb C^2$. At
the same time every quaternion $u$ corresponds to a complex duplex
matrix, i.e.,
$$
u=(a+bi)+(c+di)j\longmapsto \left( \begin{array}{cc} a+bi&c+di\\
-c+di&a-bi\end{array}\right).
$$
In \cite{Br}, $Sp(1)=\{u\in \mathbb Q: |u|=1\}$ is called the
quaternion group, or group of unit quaternions. It is known that $Sp(1)\cong SU(2)$
and is the two-fold covering group of $SO(3)$.

The automorphism of Heisenberg type groups was given by Kaplan and
Ricci \cite{Kar}. We define the translation and dilation operators
 respectively by
 $$
 \begin{array}{rl}
 T_{(x,t)}: (x^{\prime},t^{\prime})\mapsto
 (x+x^{\prime},t+t^{\prime}-2\Im(\overline {x^{\prime}}x)),
 \,\,\,\,\, (x,t),(x^{\prime},t^{\prime})\in \mathscr Q
 \end{array}
 $$
and
 $$
 \begin{array}{rl}
 T_{\rho}: (x^{\prime}, t^{\prime})\mapsto
 (\sqrt \rho x^{\prime}, \rho t^{\prime}),
 \,\,\,\,\, (x^{\prime}, t^{\prime})\in \mathscr Q, \rho>0.
 \end{array}
 $$
Let
 $$
 \begin{array}{rl}
 A_{u,v}(q)=uq\overline v,\quad u,v\in Sp(1), q \in \mathbb Q
 \end{array}
 $$
and
 $$
 \begin{array}{rl}
B_v(r)=vr\overline v,\quad v\in Sp(1), r \in \Im \mathbb Q.
\end{array}
 $$The maps
$A_{u,v}$ act transitively on the unit sphere of $\mathscr Q$, and
$B_v$ act transitively on the unit sphere of $\Im \mathbb Q$. We define the operator $T_{u,v}$ on
 $\mathscr Q$  by
 $$
 \begin{array}{rl}
 T_{u,v}: (x^{\prime},t^{\prime})\mapsto
 (ux^{\prime}\overline v, vt^{\prime}\overline v),
 \,\,\,\,\, (x^{\prime},t^{\prime})\in \mathscr Q.
 \end{array}
 $$
We are now in a position to give the affine automorporphism group of $\mathscr Q$. Let
$$
\mathbf P=\left\{(x,t,\rho,u,v): (x,t)\in \mathscr Q, \rho>0,
u,v\in Sp(1)\right\}.
$$
The action of $\mathbf P$ on $\mathscr Q$ is given by
$$
(x,t,\rho,u,v)(x^{\prime},t^{\prime})=\left(x+\sqrt \rho u
x^{\prime}\overline v, t+\rho  v t^{\prime}\overline v -2\sqrt
{\rho}\Im (v\overline{x^{\prime}}\overline {u}x)\right).
$$
That is $(x,t,\rho,u,v)= T_{(x,t)}T_{\rho}T_{u,v}$. It is known that $\mathbf P$ is two-fold covering of the affine automorporphism group of $\mathscr Q$.
The group law of $\mathbf P$ is given by
$$
\begin{array}{rl}
&(x_1,t_1,\rho_1,u_1,v_1)(x_2,t_2,\rho_2,u_2,v_2)\\[4mm]
&\quad =(x_1+\sqrt \rho_1u_1x_2\overline v_1, t_1+\rho_1
v_1t_2\overline v_1 -2\sqrt {\rho_1}\Im (v_1\overline {x_2}\,
\overline {u_1} x_1), \rho_1\rho_2, u_1u_2,v_1v_2).
\end{array}
$$
We will consider $\mathbf P$
instead of the affine automorporphism group of $\mathscr Q$.
It is easy to verify that $\mathbf P$ is a locally compact
non-unimodular group with the left Haar measure $dm_l(x,t,\rho,u,v)=\frac {dxdtd\rho dudv}{\rho^6}$ and the right Haar measure
$dm_r(x,t,\rho,u,v)=\frac {dxdtd\rho dudv}{\rho}$ respectively, where
$du$ and $dv$ are the normalized Haar measures of group $Sp(1)$.

Let us consider the unitary representation $U$ of $\mathbf P$ on $L^2(\mathscr Q)$
defined by
$$
U(x,t,\rho,u,v)f(x^{\prime},t^{\prime})=\rho^{-5/2}f\left(\frac
{{\overline u}(x^{\prime}-x){v}}{\sqrt \rho}, \frac {{\overline
v}(t^{\prime}-t+2\Im (\overline
x^{\prime}x))v}{\rho}\right).
$$
The representation $U$ is reducible on $L^2(\mathscr Q)$. We shall
decompose the space $L^2(\mathscr Q)$ into the direct sum of the
irreducible invariant closed subspaces.

\vskip 0.3 true cm
\section{\bf  Direct sum decomposition for $L^2(\mathscr Q)$}
\vskip 0.3 true cm

First we state some results of the Fourier
analysis on  $\mathscr Q$. The Fourier transform on Heisenberg type groups was studied by Kaplan and
Ricci \cite{Kar}. Let $0\not= a\in \Im \mathbb Q$. Set
$\tilde{a}=\frac {a}{|a|}$. The mapping $\rho(\tilde{a}):
q\longrightarrow q\tilde a$ gives a complex structure
of $\mathbb Q$. Let $\mathscr H_{a}$ be the Fock space consisting
of all holomorphic functions $F$ on $(\mathbb Q,
\rho(\tilde{a}))\cong \mathbb {C}^2$ such that
$$
\|F\|^2=\int_{\mathbb Q}|F(q)|^2e^{-2|a||q|^2}dq< \infty.
$$
We now define the unitary representation $\pi_{a}(x,t)$ of
$\mathscr Q$ on $\mathscr H_{\alpha}$ by
$$
\pi_{a}(x,t)F(q)=F(q+x)e^{i\langle a,t\rangle-|a|(|x|^2+2\langle
q,x\rangle -2i\langle q\tilde{a},x\rangle)}.
$$
Up to a unitary equivalent, all irreducible infinite-dimensional
unitary representations of $\mathscr Q$ are given by
$\pi_{a}(x,t)$. For $f\in L^1(\mathscr Q)$, the  Fourier transform
of $f$ is an operator valued function defined by
$$
\widehat{f}(a)=\int_{\mathscr Q}f(x,t)\pi_{a}(x,t)dxdt.
$$
Let $f,g\in L^1(\mathscr Q)\cap L^2(\mathscr Q)$. By the standard
theory of the Weyl transform,  we have
$$
\langle f,g\rangle_{L^2(\mathscr Q)}=\displaystyle{\frac
{1}{2\pi^5}\int_{\Im \mathbb Q }\mbox {tr}(\widehat g(a)^{\ast}
\widehat f(a))|a|^2da}.
$$
Specially, the following Plancherel formula holds.
$$
\|f\|^2_{L^2(\mathscr Q)}=\displaystyle{\frac {1}{2\pi^5}\int_{\Im
\mathbb Q }\|\widehat f(a)\|^2_{HS}|a|^2da.}\eqno(2.1)
$$
The Fourier transform can be extended to the tempered
distributions on $\mathscr Q$ by duality. And we also have the
formula of inverse Fourier transform
$$
f(x,t)=\displaystyle{\frac {1}{2\pi^5}\int_{\Im \mathbb Q }\mbox
{tr}(\widehat \pi_{a}^{\ast}(x,t) \widehat
f(a))|a|^2da}.
$$
Let $f*g$ denote the convolution of $f$ and $g$, i.e.,
$$
f*g(x,t)=\displaystyle{\int_{\mathscr
Q}f(y,s)g((y,s)^{-1}(x,t))dyds}.
$$
Then
$$
\widehat {f*g}(a)=\widehat f(a)\widehat g(a).\eqno(2.2)
$$
If $\tilde
f(x,t)=\overline {f((x,t)^{-1})}=\overline {f(-x,-t)}$, then
$$
\widehat {\tilde f}(a)=\widehat f(a)^{\ast},\eqno(2.3)
$$
where $\widehat f(a)^{\ast}$ is the adjoint of $\widehat f(a)$.
We now choose an orthonormal basis $\{e_0,e_1,e_2,e_3\}$ of
$\mathbb Q$ satisfying $e_0=1, e_1=\tilde{a}, e_2\tilde{a}=e_3$.
Write $q=(z_1,z_2)$ for $z_1=x_1+y_1i, z_2=x_2+y_2i$ and
$q=x_1e_0+y_1e_1+x_2e_2+y_2e_3$. Then
$$
\displaystyle{\{E_{\alpha}^a
(q)=\pi^{-1}(\alpha_1!\alpha_2!)^{-1/2}(2|a|)^{\frac
{\alpha_1+\alpha_2}{2}+1}z_1^{\alpha_1}z_2^{\alpha_2}:\alpha=(\alpha_1,\alpha_2)\in
\mathbb N^2\}}
$$
is an orthonormal basis of $\mathscr H_{a}$. Thus we have the
identification between Fock spaces $\mathscr H_{a}$ and $\mathscr
H_{a^{\prime}}$ by identifying $E_{\alpha}^a $ with
$E_{\alpha}^{a^{\prime}}$. For $\rho>0$, write
$$
f_{\rho}(x,t)=\rho^{-5}f(\frac {x}{\sqrt \rho}, \frac {t}{\rho}),
$$
then we have
$$
\widehat f_{\rho}(a)=\widehat f(\rho a),
$$
where we have identified $F(q) \in \mathscr H_{a}$ with $q F(\sqrt \rho q) \in \mathscr H_{\rho a}$.
Suppose that $u, v\in Sp(1)$ and $f_{u,v}(x,t)=f(\overline ux
v,\overline vt v)$, we can verify that
$$
\begin{array}{rl}
\widehat f_{u,v}(a)F(q) =&\displaystyle{\int_{\mathscr Q}f(\overline u x
v, \overline vt v)e^{i\langle a,t\rangle-|a|(|x|^2+2\langle q,
x\rangle-2i\langle q\widetilde
a,x\rangle)}F(q+x)dxdt}\\[4mm]
=&\displaystyle{\int_{\mathscr
Q}f(x,t)e^{i\langle \overline va v
,t \rangle-|\overline va v|(|x|^2+2\langle
\overline uq  v, x \rangle-2i\langle (\overline uq v)(
\overline v\widetilde a v), x \rangle)}F(u(\overline u q
v+x )\overline v) dx dt}.
\end{array}
$$
Thus,
$$
\widehat f_{u,v}(a)=\gamma_{u,v}^{-1} \widehat f(\overline va v
)\gamma_{u,v},
$$
where the intertwining operator $\gamma_{u,v}$ is given by
$$
\gamma_{u,v}F(q)=F(uq \overline v ).
$$
Let $l=\alpha_1+\alpha_2$ and $\mathscr H_{a,l}$ be the subspace of
$\mathscr H_a$ which consists of all homogeneous polynomials of
degree $l$ in $q\in \mathbb C^2$. Then $\mathscr H_{a,l}$ is an
irreducible invariant closed subspace under $\gamma_{u,v}$ where we
have identified $F(q) \in \mathscr H_{a}$ with $F(q \overline v) \in
\mathscr H_{\overline va v}$. Moreover, we have
$$
\mathscr H_a=\bigoplus_{l=0}^{+\infty}\mathscr H_{a,l}.
$$
Let $P_{a,l}$ denote the orthogonal projection operator from
$\mathscr H_a$ to $\mathscr H_{a,l}$. The projection operator $P_l$
is defined in terms of the Fourier transform by
$$
\widehat {P_lf}(a)=\widehat f(a)P_{a,l}.
$$
Clearly, $P_l$'s are mutual orthogonal projection operators on $L^2(\mathscr
Q)$.
If $f(x,t)$ is a
radial function with respect to the variable
$x$, i.e., $f(ux,t)=f(x, t)$ for all $u \in Sp(1)$, then we have $\gamma_{u,1}\widehat f(a)=\widehat f(a)
\gamma_{u,1}$. By Schur's lemma,
$$
\widehat f(a)=\sum_{l=0}^{\infty}B_f(a,l)P_{a,l},\eqno(2.4)
$$
where $B_f(a,l)$ is a constant depending on $f, a$ and $l$. Define the subspace $H_l$ of $L^2(\mathscr Q)$ by
$$
H_l=\displaystyle{\{f\in L^2(\mathscr Q): \widehat f(a)=\widehat
f(a) P_{a,l}\}}
$$
which is the range of projection operator $P_l$. For a function $f\in L^2(\mathscr Q)$, it is easy to verify that
$$
\displaystyle{\big({U(x,t,\rho,u,v)f}\big )\widehat{\ \ }(a)=\rho^{5/2}\pi_a(x,t)\widehat {f_{u,v}}(\rho
 a).}\eqno(2.5)
$$
By a similar argument as in \cite {JL4}, we obtain the direct sum decomposition
for $L^2(\mathscr Q)$ as follows.

{\bf Theorem 1.} {\it $H_l$ is an irreducible invariant
closed subspace of $L^2(\mathscr Q)$ under the unitary representation $U$ of $\mathbf P$, and we have
$$
L^2(\mathscr
Q)=\displaystyle{\bigoplus_{l=0}^{\infty}H_l}.
$$}

{\bf Remark 1.} {\it In fact, $H_l$ can be characterized by the
subLaplacian operator $\Delta_{\mathscr Q}$ (see \cite{Tie}), i.e.,
$f \in H_l$ if and only if
$$
\widehat {\Delta_{\mathscr Q}f}(a)=-8(l+1)|a|\widehat
f(a),
$$
where $\Delta_{\mathscr Q}=X_0^2+X_1^2+X_2^2+X_3^2$ is the square sum of horizontal vector fields, $X_0, X_1,
X_2$ and $X_3$ are left invariant vector fields given by
$$
\begin{array}{rl}
&X_0=\displaystyle{\frac {\partial}{\partial x_0}-2x_1\frac
{\partial}{\partial t_1}-2x_2\frac {\partial}{\partial
t_2}-2x_3\frac {\partial}{\partial t_3}}\\[4mm]
&X_1=\displaystyle{\frac {\partial}{\partial x_1}+2x_0\frac
{\partial}{\partial t_1}-2x_3\frac {\partial}{\partial
t_2}+2x_2\frac {\partial}{\partial t_3}}\\[4mm]
&X_2=\displaystyle{\frac {\partial}{\partial x_2}+2x_3\frac
{\partial}{\partial t_1}+2x_0\frac {\partial}{\partial
t_2}-2x_1\frac {\partial}{\partial t_3}}\\[4mm]
&X_3=\displaystyle{\frac {\partial}{\partial x_3}-2x_2\frac
{\partial}{\partial t_1}+2x_1\frac {\partial}{\partial
t_2}+2x_0\frac {\partial}{\partial t_3}}.
\end{array}
$$}

\vskip 0.3 true cm
\section{\bf Continuous wavelet transforms }
\vskip 0.3 true cm

We are going to show that the restriction of $U$ on $H_l$ is
square-integrable. In other words, there exists a non-zero
function $\phi\in H_l$, such that
$$
C_{\phi}=\displaystyle{\frac {1}{\|\phi\|^2_{L^2(\mathscr
Q)}}\int_{\mathbf P} |\langle
\phi,U(x,t,\rho,u,v)\phi\rangle_{L^2(\mathscr Q)}|^2
dm_l(x,t,\rho,u,v)< \infty}.\eqno(3.1)
$$
We call (3.1) the admissibility condition, and write $\phi\in AW_l$.

{\bf Theorem 2.} {\it Let $\phi\in H_l$, not identically zero. Then
$\phi\in AW_l$ if and only if
$$
C_{\phi}=\displaystyle{\frac {1}{ d_l}\int_{\Im \mathbb Q}\|\widehat
\phi(a)\|^2_{HS}\frac {da}{|a|^3} < \infty.}\eqno(3.2)
$$
where $d_l=\text{dim}{\mathscr H_{a,l}}=l+1$. $(3.2)$ is called the
Calder\'on reproducing condition.}

{\bf Proof}. Suppose $\phi\in H_l$, by (2.5), we have
$$
\displaystyle{\int_{\mathscr Q}\langle
\phi,U(x,t,\rho,u,v)\phi\rangle_{L^2(\mathscr
Q)}\pi_a(x,t)dxdt=\rho^{5/2}\widehat \phi(a)\widehat
{\phi_{u,v}}(\rho a)^*}.
$$
Using the Plancherel formula (2.1) we can derive
$$
\begin{array}{rl}
&\displaystyle{\int_{\mathbf P}|\langle
\phi,U(x,t,\rho,u,v)\phi\rangle_{L^2(\mathscr Q)}|^2dm_l(x,t,\rho,u,v)}\\[4mm]
&\quad =\displaystyle{\frac
{1}{2\pi^5}\int_{Sp(1)}\int_{Sp(1)}\int_{0}^{+\infty}\left(\int_{\Im
\mathbb Q}\|\widehat \phi(a)\widehat {\phi_{u,v}}(\rho
a)\|_{HS}^2|a|^2da\right)\frac {d\rho du dv}{\rho}}\\[4mm]
&\quad =\displaystyle{\frac
{1}{2\pi^5}\int_{Sp(1)}\int_{Sp(1)}\int_{0}^{+\infty}\left(\int_{\Im
\mathbb Q}\text{tr}\left(\widehat \phi(a)^*\widehat \phi(a)\widehat
{\phi_{u,v}}(\rho a)^*\widehat {\phi_{u,v}}(\rho
a)\right)|a|^2da\right)\frac {d\rho du dv}{\rho}}.
\end{array}
$$
Because
$$
\begin{array}{rl}
&\displaystyle{\text{tr}\left(\widehat \phi(a)^*\widehat
\phi(a)\widehat {\phi_{u,v}}(\rho a)^*\widehat {\phi_{u,v}}(\rho
a)\right)}\\[4mm]
&\quad =\displaystyle{\sum_{|\alpha|=l}\left\langle \widehat
{\phi_{u,v}}(\rho a)^*\widehat {\phi_{u,v}}(\rho a)E_{\alpha}^a ,
\widehat \phi(a)^*\widehat \phi(a)E_{\alpha}^a
\right\rangle_{\mathscr H_a}},
\end{array}
$$
 we therefore obtain
$$
\begin{array}{rl}
&\displaystyle{\int_{\mathbf P}|\langle
\phi,U(x,t,\rho,u,v)\phi\rangle_{L^2(\mathscr Q)}|^2dm_l(x,t,\rho,u,v)}\\[4mm]
&\quad =\displaystyle{\frac
{1}{2\pi^5}\int_{0}^{+\infty}\left\{\int_{\Im \mathbb
Q}\sum_{|\alpha|=l}\left\langle\left(\int_{Sp(1)\times
Sp(1)}\widehat {\phi_{u,v}}(\rho a)^*\widehat {\phi_{u,v}}(\rho
a)dudv\right)E_{\alpha}^a ,\right.\right.}\\[4mm]
&\hskip 7 true cm \displaystyle{\left.\left. \widehat
\phi(a)^*\widehat \phi(a)E_{\alpha}^a \right\rangle_{\mathscr
H_a}|a|^2dv da\right\}\frac {d\rho}{\rho}.}
\end{array}
$$
 By (2.2) together with (2.3), we derive
$$
\displaystyle{\int_{Sp(1)\times Sp(1)}\widehat {\phi_{u,v}}(\rho
a)^*\widehat {\phi_{u,v}}(\rho
a)dudv}=\displaystyle{\int_{Sp(1)\times Sp(1)}\left(\widehat
{\widetilde {\phi}_{u,v}*\phi_{u,v}}\right)(\rho
a)dudv}=\displaystyle{\widehat {\psi}(\rho a)},
$$
where
$$
\begin{array}{rl}
\psi(x,t)&=\displaystyle{\int_{Sp(1)\times Sp(1)}\left({\widetilde
{\phi}_{u,v}*\phi_{u,v}}\right)(x,t)dudv}\\[4mm]
&=\displaystyle{\int_{Sp(1)\times Sp(1) }\left({\widetilde
{\phi}*\phi}\right)(\overline ux v, \overline vt
v)dudv}\\[4mm]
&=\displaystyle{\int_{Sp(1)}\vartheta_v( x, t)dv}
\end{array}
$$
$$
\vartheta(x,t)=\int_{Sp(1)}{\widetilde {\phi}*\phi}(\overline
ux,t)du    \,\, \text{and}\,\, \vartheta_v(x,t)=\vartheta(x
v,\overline v t v).
$$
 Obviously, $\vartheta$ is a radial function with respect to the variable
$x$. Hence by (2.4) we get $\widehat \vartheta(\rho
a)=B_{\vartheta}(\rho a)P_l$. The Fourier transform of $\psi$ is
given by
$$
\widehat \psi(\rho a)=\displaystyle{\int_{Sp(1)}B_{\vartheta}
(\rho\overline v a v)P_ldv}.
$$
Consequently, we have
$$
\begin{array}{rl}
&\displaystyle{\int_{\mathbf P}|\langle
\phi,U(x,t,\rho,u,v)\phi\rangle_{L^2(\mathscr Q)}|^2dm_l(x,t,\rho,u,v)}\\[4mm]
&\quad =\displaystyle{\frac
{1}{2\pi^5}\int_{0}^{+\infty}\left\{\int_{\Im \mathbb
Q}\int_{Sp(1)}\sum_{|\alpha|=l}\left\langle\left(B_{\vartheta}(\rho\overline
v av) \widehat \phi(a)^*\widehat \phi(
a)\right)E_{\alpha}^a ,E_{\alpha}^a \right\rangle_{\mathscr H_a}|a|^2da\right\}\frac {dv d\rho}{\rho}}\\[4mm]
&\quad =\displaystyle{\frac {1}{2\pi^5}\int_{0}^{+\infty}\int_{\Im
\mathbb Q}\int_{Sp(1)}B_{\vartheta}(\rho\overline v
av)\text{tr}\left(\widehat \phi(a)^*\widehat \phi(
a)\right)|a|^2da\frac {dv
d\rho}{\rho}}\\[4mm]
&\quad
=\displaystyle{\left(\int_{Sp(1)}\int_{0}^{+\infty}B_{\vartheta}
(\rho\overline v av)\frac {dv d\rho}{\rho}\right)\left(\frac
{1}{2\pi^5}\int_{\Im \mathbb Q}\text{tr}(\widehat \phi(a)^*\widehat
\phi(a))|a|^2da\right)}\\[4mm]
&\quad =\displaystyle{\|\phi\|^2_{L^2(\mathscr Q)}\left(\int_{\Im
\mathbb Q}B_{\vartheta}(a)\frac {da}{|a|^3}\right)}.
\end{array}
$$
Because
$$
\begin{array}{rl}
B_{\vartheta}( a)&=\displaystyle{\frac {1}{d_l}\text{tr}\widehat
\vartheta(a)}\\[4mm]
&=\displaystyle{\frac {1}{d_l}\text{tr} \left( \Big(\int_{Sp(1)}
{\phi}* {\phi}(\overline u x,t)du\Big)\widehat{\ \ }(a)\right)}\\[4mm]
&=\displaystyle{\frac {1}{d_l}\text{tr}\big(\widehat
{\phi}(a)^*\widehat {\phi}(a)\big)}\\[4mm]
&=\displaystyle{\frac {1}{d_l}\|\widehat \phi(a)\|_{HS}^2},
\end{array}
$$
Theorem 2 is proved. \qed

Let $\phi\in AW_l, f\in H_l$. The continuous wavelet transform
$W_{\phi}$ on $H_l$ is defined by
$$
W_{\phi}f(x,t,\rho,u,v)=\langle f,
U(x,t,\rho,u,v)\phi\rangle_{L^2(\mathscr Q)}.
$$
For $\phi, \psi \in AW_l$, we set
$$
\langle \phi, \psi\rangle_{AW_l} =\displaystyle{\frac {1}{
d_l}\int_{0}^{+\infty}\mbox{tr}(\widehat \psi( a)^*\widehat \phi(
a))\frac {da}{|a|^3}}
$$
and call it the inner product in $AW_l$.

{\bf Theorem 3.} {\it
Let $\phi\in AW_l, \psi \in AW_{l^{\prime}}$, $f\in H_l, g\in
H_{l^{\prime}}$. Then
$$
\langle W_{\psi}f, W_{\phi}g\rangle_{L^2(\mathbf P, dm_l)}=\langle
\phi,\psi\rangle_{AW_l}\langle f,g\rangle_{L^2(\mathscr Q)}.
$$
Specially,
$$
\|W_{\phi}f\|_{L^2(\mathbf
P,dm_l)}=C_{\phi}^{1/2}\|f\|_{L^2(\mathscr Q)}
$$
and
$$
\langle W_{\psi}f, W_{\phi}g\rangle_{L^2(\mathbf P, dm_l)}=0
$$
when $l\not=l^{\prime}$.}

{\bf Proof.} Let $\phi\in AW_l$, $f\in H_l$. Then
$$
\int_{\mathscr Q}\langle f,
U(x,t,\rho,u,v)\phi\rangle_{L^2(\mathscr
Q)}\pi_a(x,t)dxdt=\rho^{5/2}\widehat f(a)\widehat
{\phi_{u,v}}(\rho a)^{\ast}.
$$
By the Plancherel formula we can get
$$
\begin{array}{rl}
&\displaystyle{\bigg\langle  W_{\psi}f,W_{\phi}g \bigg\rangle_{L^2(\mathbf P,dm_l)}}\\[4mm]
&\quad =\displaystyle{\int_{\mathbf P}W_{\psi}f(x,t,\rho,u,v)
\overline {W_{\phi}g(x,t,\rho,u,v)}dm_l(x,t,\rho,u,v)}\\[4mm]
&\quad =\displaystyle{\frac
{1}{2\pi^5}\int_{Sp(1)}\int_{Sp(1)}\int_0^{+\infty}\left(\int_{\Im
\mathbb Q}\mbox{tr}(\widehat g(a)^{\ast}\widehat f(a)\widehat
{\psi_{u,v}}(\rho a)^{\ast}\widehat \phi_{u,v}(\rho
a))|a|^2da\right)\frac {d\rho
du dv}{\rho}}\\[4mm]
&\quad =\displaystyle{\frac
{1}{2\pi^5}\int_0^{+\infty}\left\{\int_{\Im \mathbb
Q}\sum_{\alpha\in \mathbb N^2}\bigg\langle\left(\int_{Sp(1)\times
Sp(1)}\widehat \psi_{u,v}(\rho a)^{\ast}\widehat \phi_{u,v}(\rho
a)dudv\right)E_{\alpha}^a ,\right.}\\[4mm]
&\displaystyle{\hskip 7 true cm \left. \widehat g(a)^{\ast}\widehat
g(a)^{\ast}\widehat f(a)E_{\alpha}^a \bigg\rangle_{\mathscr
H_{a}}|a|^2da\right\}\frac { d\rho}{\rho}}.
\end{array}
$$
Suppose $\psi,\phi\in AW_l, f,g\in H_l$. Then
$$
\displaystyle{\int_{Sp(1)\times Sp(1)}\widehat \psi_{u,v}(\rho
a)^{\ast} \widehat \phi_{u,v}(\rho a)dudv}=\displaystyle{\frac
{1}{d_l}\text{tr}(\widehat{\psi}(\rho a)^*\widehat \phi(\rho
a))P_l.}
$$
Hence we get
$$
\begin{array}{rl}
&\displaystyle{\bigg\langle  W_{\psi}f,W_{\phi}g \bigg\rangle_{L^2(\mathbf P,dm_l)}}\\[4mm]
&\quad =\displaystyle{\left(\frac 1{d_l}\int_{\Im \mathbb Q}
\text{tr}(\widehat \psi(a^{\prime})^{\ast}\widehat \phi( a^{\prime})
)\frac
{da^{\prime}}{|a^{\prime}|^3}\right)\left(\frac{1}{2\pi^5}\int_{\Im
\mathbb Q}\text{tr}(\widehat g(a)^*\widehat f(a))|a|^2da\right)}\\[4mm]
&\quad =\displaystyle{\langle \phi,\psi\rangle_{AW_l}\langle
f,g\rangle_{L^2(\mathscr Q)}}.
\end{array}
$$
If $\psi\in AW_l, \phi\in AW_{l^{\prime}}, f\in
H_l, g\in H_{l^{\prime}}, l\not= l^{\prime}$, then
$$
\displaystyle{\int_{Sp(1)\times Sp(1)}\widehat \psi_{u,v}(\rho
a)^{\ast} \widehat \phi_{u,v}(\rho a)dudv}=0.
$$
Theorem 3 is proved. \qed

Theorem 2 and Theorem 3 restrict wavelets to subspaces $H_l$. The
restriction is removable. Suppose $\phi \in L^2(\mathscr{Q})$. By
Theorem 1,
$$\phi= \sum_{l=0}^{\infty} \phi_l,$$
where $\phi_l \in H_l$. If there exists a constant $C_\phi$, which
is independent of $l$, such that
$$
\frac {1}{ d_l}\int_{\Im \mathbb Q}\|\widehat
\phi_l(a)\|^2_{HS}\frac {da}{|a|^3}= C_\phi < \infty \quad \text{for
all}\,\, l\in \mathbb{N},
$$
we say that $\phi$ is an admissible wavelet and write $\phi \in AW$.
Then we define the continuous wavelet transform of $f \in
L^2(\mathscr Q)$ by
$$
W_{\phi}f(x,t,\rho,u,v)=\langle f,
U(x,t,\rho,u,v)\phi\rangle_{L^2(\mathscr Q)}.
$$

{\bf Theorem 4.} {\it Suppose $\phi \in AW$. Then
$$
\|W_{\phi}f\|_{L^2(\mathbf
P,dm_l)}=C_{\phi}^{1/2}\|f\|_{L^2(\mathscr Q)}.\eqno(3.3)
$$
}

{\bf Proof.} By Theorem 1,
$$f= \sum_{l=0}^{\infty} f_l,$$
where $f_l \in H_l$. By Theorem 3,
$$
\begin{array}{rl}
\displaystyle{\|W_{\phi}f\|_{L^2(\mathbf P,dm_l)}^2}&=\displaystyle{
\sum_{l=0}^{\infty} \|W_{\phi_l}f_l\|_{L^2(\mathbf P,dm_l)}^2}\\[4mm]
&=\displaystyle{\sum_{l=0}^{\infty} C_\phi \|f_l\|_{L^2(\mathscr Q)}^2}\\[4mm]
&=\displaystyle{C_\phi \|f\|_{L^2(\mathscr Q)}^2.}\\[4mm]
\end{array}
$$
This proves Theorem 4. \qed

The fundamental manifold of the quaternion Heisenberg group
$\mathscr Q$ is just $\mathbb R^4\times \mathbb R^3$. The Schwartz
space $\mathscr S(\mathscr Q)$ coincides with the Schwartz space on
$\mathbb R^4\times \mathbb R^3$. As a consequence of Theorem 4, the
inverse wavelet transform is valid. We state the result as follows.

{\bf Theorem 5.} {\it Suppose $\phi \in AW$ and $f \in L^2(\mathscr
Q)$. Then
 $$
 \displaystyle{f(x,t)=C_{\phi}^{-1}\int_{\mathbf P}
W_{\phi}f(y,s,\rho,u,v))U(y,s,\rho,u,v)\phi(x,t)dm_l(y,s,\rho,u,v)}
$$
in $L^2$-sense. Furthermore, if $f\in \mathscr S(\mathscr Q)$, then
the above formula holds pointwise.}

The proof of Theorem 5 is standard.

Now we construct a class of admissible wavelets. These admissible wavelets are radial.
The inverse wavelet transform can be simplified by using radial wavelets.

Let $\eta$ be a function on
$\mathbb{R}_+$ satisfying
$$\int_0^{\infty} |\eta(r)|^2 r^4 dr < \infty$$
and the function $\phi_\eta$ is defined by
$$\widehat \phi_\eta (a)= \sum_{l=0}^{\infty} \eta\big( (l+1)|a|\big)
P_{a,l}.$$ Then $\phi_\eta \in L^2(\mathscr{Q})$ because, by the
Plancherel formula,
$$
\begin{array}{rl}
\displaystyle{\|\phi_\eta\|^2_{L^2(\mathscr Q)}}
&=\displaystyle{\frac {1}{2\pi^5}\int_{\Im \mathbb Q }\|\widehat \phi_\eta(a)\|^2_{HS}|a|^2da}\\[4mm]
&=\displaystyle{\frac {1}{2\pi^5}\int_{\Im \mathbb Q }
\left(\sum_{l=0}^{\infty} \big|\eta\big( (l+1)|a|\big)\big|^2 (l+1)\right)|a|^2da}\\[4mm]
&=\displaystyle{\frac {1}{2\pi^5} \sum_{l=0}^{\infty}
\frac{1}{(l+1)^4}\int_{\Im \mathbb Q } |\eta(|a|)|^2 |a|^2da}\\[4mm]
&=\displaystyle{C \int_0^{\infty} |\eta(r)|^2 r^4 dr < \infty.}\\[4mm]
\end{array}
$$
Furthermore, if $\eta$ satisfies the Calder\'on reproducing
condition
$$C_\eta= 4 \pi \int_0^{\infty} |\eta(r)|^2 \frac{dr}{r} < \infty, \eqno(3.4)$$
then $\phi_\eta$ is an admissible wavelet. In fact, write
$$\phi_\eta= \sum_{l=0}^{\infty} \phi_{\eta,l}$$
where $\phi_{\eta,l} \in H_l$, it is clear that
$$\widehat \phi_{\eta,l}(a)= \eta\big( (l+1)|a|\big) P_{a,l}.$$
Then, for all $l \in \mathbb{N}$,
$$
\begin{array}{rl}
\displaystyle{\frac {1}{ d_l}\int_{\Im \mathbb Q}\|\widehat
\phi_{\eta,l}(a)\|^2_{HS}\frac {da}{|a|^3}}
&=\displaystyle{\int_{\Im \mathbb Q }
\big|\eta\big( (l+1)|a|\big)\big|^2 \frac {da}{|a|^3}}\\[4mm]
&=\displaystyle{\int_{\Im \mathbb Q } |\eta(|a|)|^2 \frac {da}{|a|^3}}\\[4mm]
&=\displaystyle{C_\eta < \infty.}\\[4mm]
\end{array}
$$
Note that the admissible wavelet $\phi_\eta$ defined above is a radial
function, i.e., $(\phi_\eta)_{u,v}=\phi_\eta$. The wavelet transform
$W_{\phi_\eta}f$ is independent of $u$ and $v$, i.e.,
$$
W_{\phi_\eta}f(x,t, \rho, u,v)=W_{\phi_\eta}f(x,t, \rho, 1,1).
$$
Thus the integration over $u$ and $v$ can be reduced if we use the admissible wavelet $\phi_\eta$.
Write $W_{\phi_\eta}f(x,t, \rho)$ and $U(x,t,\rho)$ instead of $W_{\phi_\eta}f(x,t, \rho, 1,1)$ and $U(x,t,\rho,1,1)$.
Then we have

{\bf Theorem 6.} {\it Suppose the admissible wavelet $\phi_\eta$ is defined as above, and $f \in L^2(\mathscr
Q)$. Then
$$
\int_{\mathscr Q\times \mathbb R_+} |W_{\phi_\eta}f (x,t, \rho)|^2 \frac
{dxdtd\rho}{\rho^6} =C_{\eta}^{1/2}\|f\|_{L^2(\mathscr Q)},
$$
and
$$
 \displaystyle{f(x,t)=C_{\eta}^{-1}\int_{\mathscr Q\times \mathbb R_+}
W_{\phi_\eta}f(y,s,\rho))U(y,s,\rho)\phi_\eta(x,t) \frac
{dydsd\rho}{\rho^6}} \eqno(3.5)
$$
in $L^2$-sense. Furthermore, if $f\in \mathscr S(\mathscr Q)$, then
$(3.5)$ holds pointwise.}

\vskip 0.3 true cm
\section{\bf The Radon transform }
\vskip 0.3 true cm

Similar to the Radon transform on the Heisenberg group $\mathbf H_n$
defined by Strichartz \cite{Str}, we define the Radon transform
$\mathscr R$ on $\mathscr Q$ by
$$
\begin{array}{rl}
\mathscr R(f)(x,t)&=\displaystyle{\int_{\mathbb Q}f((x,t)(y,0))dy}\\[4mm]
&=\displaystyle{\int_{\mathbb Q}f(y,t-2\Im (\overline y x))dy}.
\end{array}
$$
Even if $f\in \mathscr S(\mathscr Q)$, the Radon transform $\mathscr
R(f)$ may not be rapidly deceasing at infinity. This fact were
pointed by Strichartz \cite{Str} in the case of $\mathbf H_n$. We
will find a subspace of Schwartz space $\mathscr S(\mathscr Q)$, on
which the Radon transform is a bijection.

Let $\mathscr F_2$ denote the Euclidean Fourier transform with
respect to the central variable $t$ alone and $\mathscr F$ denote
the full Euclidean Fourier transform, i.e.,
$$\begin{array}{rl}
\mathscr F_2(f)(x,a)&=\displaystyle{\int_{\Im \mathbb Q}f(x,t)e^{i\langle
t, a\rangle}dt,}\\[4mm]
\mathscr F(f)(y,a)&=\displaystyle{\int_{\mathscr Q}f(x,t)e^{i\langle
x,y\rangle}e^{i\langle t,a\rangle}dxdt}.
\end{array}
$$
Because
$$
\begin{array}{rl}
\mathscr F_2(\mathscr R(f))(x,a)&=\displaystyle{\int_{\Im
\mathbb Q}\left(\int_{\mathbb Q}f(y,t-2\Im(\overline y
x))dy\right)e^{i\langle t,a\rangle}dt}\\[4mm]
&=\displaystyle{\int_{\mathbb Q}\mathscr F_2(f)(y,a)e^{2i\langle
\Im(\overline y x), a\rangle}dy}\\[4mm]
&=\displaystyle{\int_{\mathbb Q}\mathscr F_2(f)(y,a)e^{-2i\langle
y, xa\rangle}dy},
\end{array}
$$
we have
$$
\mathscr F_2(\mathscr R(f))(x,a)=\displaystyle{\mathscr
F(f)(-2xa, a)}.\eqno(4.1)
$$

Let $t^s=t_1^{s_1}t_2^{s_2}t_3^{s_3}$ where $s=(s_1,s_2,s_3)\in \mathbb
N^3$. We define the space
$\mathscr S_*(\mathscr Q)$ by
$$
\mathscr S_*(\mathscr Q)=\displaystyle{\left\{f(x,t)\in \mathscr
S(\mathscr Q): \int_{\Im \mathbb Q}f(x,t)t^sdt=0 \,\, \text{for
all}\,\, x\in \mathbb Q, s\in \mathbb N^3\right\}}.
$$
Write
$\partial^{s}_t f(x,0)=\frac {\partial^{|s|}}{\partial
{t_1^{s_1}}\partial {t_2^{s_2}}
\partial {t_3^{s_3}}}f(x,t)|_{t=0}$ and define the space
$\mathscr S^*(\mathscr Q)$ by
$$
\mathscr S^*(\mathscr Q)=\displaystyle{\left\{ f\in \mathscr
S(\mathscr Q): \partial ^s_t f(x,0)=0 \,\, \text{for
all}\,\, x\in \mathbb Q, s\in \mathbb N^3\right\}}.
$$
By same argument as Proposition 5.1 in \cite{Ru1}, $f\in
\mathscr S_*(\mathscr Q)$ if and only if $\mathscr
F(f)\in \mathscr S^*(\mathscr Q)$ and $\mathscr F_2$ is an isomorphism from
$\mathscr S_*(\mathscr Q)$ onto $\mathscr S^*(\mathscr Q)$.
The space $\mathscr S^*(\mathscr Q)$ and $\mathscr S_*(\mathscr Q)$ are regarded
as Semyanistyi-Lizonkin type spaces. This kind of spaces has many applications (see \cite{Ru2}, \cite {Sa}).

We introduce a ``mixing" map $\Psi$ which given by
$$
\Psi(f)(x,t)= f(-2xt, t).
$$
Let $f \in \mathscr S^*(\mathscr Q)$. Similar to the inequality $(5.5)$ in \cite{Ru1},
for every $p,q \in \mathbb{N}$, there is a constant $c_{p,q}$ such that
$$
|f(x,t)| \leq c_{p,q} |t|^{2p}(1+|x|^2)^{-q}.
$$
It follows that $\Psi$ is a bijection on $\mathscr S^*(\mathscr Q)$.
The inverse of $\Psi$ is given by
$$
\Psi^{-1}(f)(x,t)=\left\{\begin{array}{ll}
f \big(\frac
{xt}{2|t|^2},t \big),&\quad \mbox{for $t\not=0$}\\
0, &\quad \mbox{for $t=0$}\end{array}\right.
$$
Now $(4.1)$ reads as
$$
 \mathscr F_2(\mathscr R(f))(x,a)=\Psi (\mathscr F(f))(x,a).
$$
Thus the Radon transform $\mathscr R= \mathscr F_2^{-1} \Psi \mathscr F$ is a bijection on $\mathscr S_*(\mathscr Q)$, and we have an inversion
formula of the Radon transform as follows.

{\bf Theorem 7.} {\it Let $f\in \mathscr S_*(\mathscr Q)$. Then
$$
\mathscr R^{-1}(f)=\mathscr F^{-1}\Psi^{-1}\mathscr
F_2(f).
$$}

We can give another inversion formula of $\mathscr R$ by using of
the Fourier transform on $\mathscr Q$. First, we have
$$
\begin{array}{rl}
& (\widehat {\mathscr R(f)}(a)E_{\alpha}^a )(q)\\[4mm]
&\quad =\displaystyle{\int_{\mathscr Q}\mathscr R(f)(x,t)e^{i\langle
a,t\rangle-|a|(|x|^2+2\langle q,x\rangle-2i\langle q\tilde
a,x\rangle)}E_{\alpha}^a (q+x)dxdt}\\[4mm]
&\quad =\displaystyle{\int_{\mathscr Q}\left(\int_{\mathbb
Q}f(y,t-2\Im(\overline yx))dy\right)e^{i\langle
a,t\rangle-|a|(|x|^2+2\langle q,x\rangle-2i\langle q\tilde
a,x\rangle)}E_{\alpha}^a (q+x)dxdt}\\[4mm]
&\quad =\displaystyle{\int_{\mathbb Q \times \mathbb Q}\mathscr
F_2(f)(y,a) E_{\alpha}^a (q+x)e^{-|a|(|x|^2+2\langle
q,x\rangle-2i\langle (y+q)\tilde a,x\rangle)}dxdy}\\[4mm]
&\quad =\displaystyle{\int_{\mathbb Q}\mathscr F_2(f)(y,a)
e^{|a|(|q|^2+2i\langle q\tilde a,y\rangle)} \left (\int_{\mathbb Q}
E_{\alpha}^a (x)e^{-|a|(|x|^2-2i\langle (y+q)\tilde a,x\rangle)} dx
\right )dy}.
\end{array}
$$
Using the identity
$$
\int_{\mathbb C^2} z^{\alpha} e^{-s|z|^2} e^{i\langle z,w \rangle}
dz= \Big(\frac{\pi}{s}\Big)^2 \Big(\frac{iw}{2s}\Big)^{\alpha}
e^{-\frac{|w|^2}{4s}},
$$
and note that $E_{\alpha}^a $ is a holomorphic monomial with respect
to the complex structure $\rho(\tilde{a})$, we obtain
$$
\int_{\mathbb Q} E_{\alpha}^a (x)e^{-|a|(|x|^2-2i\langle (y+q)\tilde
a,x\rangle)} dx = (-1)^{\alpha_1+\alpha_2}
\Big(\frac{\pi}{|a|}\Big)^2 E_{\alpha}^a (q+y) e^{-|a||q+y|^2}.
$$
Then
$$
\begin{array}{rl}
& (\widehat {\mathscr R(f)}(a)E_{\alpha}^a )(q)\\[4mm]
&\quad =\displaystyle{(-1)^{\alpha_1+\alpha_2}\Big(\frac{\pi}{|a|}\Big)^2
\int_{\mathbb Q}\mathscr F_2(f)(y,a)e^{-|a|(|y|^2+2\langle
q,y\rangle-2i\langle q\tilde a,y\rangle)}E_{\alpha}^a (q+y)dy}\\[4mm]
&\quad
=\displaystyle{(-1)^{\alpha_1+\alpha_2}\Big(\frac{\pi}{|a|}\Big)^2(\widehat
f(a)E_{\alpha}^a )(q)}.
\end{array}
$$
Thus we get
$$
\widehat {\mathscr R(f)}(a)=\Big(\frac{\pi}{|a|}\Big)^2 \widehat
f(a)S,\quad f\in \mathscr S_*(\mathscr Q),\eqno(4.2)
$$
where $S=\displaystyle{\sum_{l=0}^{\infty} (-1)^l \mathcal
P_{a,l}}$.

Let $\mathscr L=\displaystyle{-\frac{1}{\pi^2} \sum_{l=1}^3 \frac
{\partial^2}{\partial t_l^2}}$. Essentially, $\mathscr L$ is the
Laplacion with respect to the central variable $t$. It is easy to
see that
$$
\displaystyle{\widehat {\mathscr L(f)}(a)=\Big(\frac{|a|}{\pi}\Big)^2 \widehat f(a)}.\eqno(4.3)
$$

Let $f\in \mathscr S(\mathscr Q)$. For any $m \in \mathbb{N}$, $\mathscr L^m(f) \in L^2(\mathscr Q)$. By $(4.3)$ and
the Plancherel formula $(2.1)$,
$$
\int_{\Im \mathbb Q}\|\widehat f(a)\|^2_{HS}|a|^{2+k}da< \infty \quad \mbox{for all}\,\, k \in \mathbb{N}.
$$
Now we define the space $\mathscr S_{\mathscr R}(\mathscr Q)$ by
$$
\displaystyle{\mathscr S_{\mathscr R}(\mathscr Q)=\left\{ f\in
\mathscr S(\mathscr Q): \int_{\Im \mathbb Q}\|\widehat
f(a)\|^2_{HS}|a|^{2-k}da<+\infty \,\,\mbox{for all}\,\, k\in
\mathbb N\right\}}.
$$
It is easy to see that the Radon transform $\mathscr R$ is a
bijection on $\mathscr S_{\mathscr R}(\mathscr Q)$. Moreover, by $(4.2)$ and $(4.3)$, we have

{\bf Theorem 8.}  {\it Let $f\in \mathscr S_{\mathscr R}(\mathscr Q)$. Then
$$\begin{array}{rl}
\| f\|_2&=\| \mathscr L \mathscr R(f)\|_2,\\[4mm]
\mathscr R^{-1}(f)&=\displaystyle{\mathscr L\mathscr R \mathscr L(f)}.
\end{array}$$}

{\bf Remark 1.} {\it We note that $\mathscr L$ is a positive operator. For $\mu \in \mathbb{R}$, $\mathscr L^{\mu}$ can be defined by
$$
\displaystyle{\widehat {\mathscr L^{\mu}(f)}(a)=\Big(\frac{|a|}{\pi}\Big)^{2\mu} \widehat f(a)}.
$$
We also have, for $f\in \mathscr S_{\mathscr R}(\mathscr Q)$,
$$\mathscr R^{-1}(f)=\mathscr L^{\mu}\mathscr R \mathscr L^{\nu}(f), \quad \mu + \nu=2.
$$}

The next theorem means that our two treatments about the Radon transform are essentially equivalent.

{\bf Theorem 9.} {\it
 $$\mathscr S_{\mathscr R}(\mathscr Q)=\mathscr S_{*}(\mathscr Q).$$}

{\bf Proof.} Note that $\widehat{f}(a)$ is just the
Weyl transform of $f^{a}= \mathscr F_2f(\cdot,a)$. By the standard theory of the
Weyl transform,
$$
\bigg( \frac{2|a|}{\pi}
\bigg)^2\|\widehat{f}(a)\|_{HS}^2=\|f^{a}\|^2_2.
$$
If $f\in \mathscr S_{\mathscr R}(\mathscr Q)$, then for all $k\in \mathbb N$,
$$
\displaystyle{ \lim_{|a| \rightarrow 0} |a|^{-k} \int_{\mathbb Q}|
\mathscr F_2f(x, a)|^2dx=0}.
$$
It follows that
$$
\displaystyle{ \int_{ \mathbb Q}| \mathscr F_2f(x, 0)|^2dx=
\lim_{|a| \rightarrow 0} \int_{\mathbb Q}| \mathscr F_2f(x,
a)|^2dx=0}.
$$
Furthermore, we have
$$
\displaystyle{ \int_{\mathbb Q} |\partial_{a_1}\mathscr F_2f(x,
0)|^2dx= \lim_{a_1 \rightarrow 0} a_1^{-2}\int_{\mathbb Q} |\mathscr
F_2f(x, (a_1,0,0))|^2dx=0}.
$$
Inductively, we obtain, for all $s \in \mathbb N^3$,
$$
\displaystyle{ \int_{\mathbb Q}| \partial_{a}^s\mathscr F_2f(x,
0)|^2dx=0}.
$$
This means $\mathscr F_2f\in \mathscr S^{*}(\mathscr Q)$, and
equivalently $f\in \mathscr S_{*}(\mathscr Q)$.

On the other hand, suppose $f\in \mathscr S_{*}(\mathscr Q)$, i.e.,
$\mathscr F_2f\in \mathscr S^{*}(\mathscr Q)$. Then
$$
\displaystyle{ \lim_{|a| \rightarrow 0} \int_{\mathbb Q} |\mathscr
F_2f(x, a)|^2dx= \int_{\mathbb Q} |\mathscr F_2f(x, 0)|^2dx=0}.
$$
And
$$
\lim_{|a| \rightarrow 0} |a|^{-2}\int_{\mathbb Q} |\mathscr F_2f(x,
a)|^2dx =\lim_{|a| \rightarrow 0} \int_{\mathbb Q} \bigg|
\sum_{l=1}^3 \frac{a_l}{|a|}
\partial_{a_l} \mathscr F_2f(x, \theta a)\bigg|^2dx=0,
$$
where $0< \theta <1$. Inductively, for all $k\in \mathbb N$,
$$
\displaystyle{ \lim_{|a| \rightarrow 0} |a|^{-k} \int_{\mathbb Q}|
\mathscr F_2f(x, a)|^2dx=0}.
$$
Therefore, $f\in \mathscr S_{\mathscr R}(\mathscr Q)$. Theorem 9 is proved. \qed

\vskip 0.3 true cm
 \section  {\bf Inverse Radon transform by using
wavelets} \vskip 0.3 true cm

The inversion formulas of the Radon transform in above section require the smoothness of functions.
In this section we establish an inversion formula of the Radon transform in $L^2$-sense by using the inverse wavelet
transform. The formula does not require the smoothness of functions. Instead, we will use smooth wavelets.

Suppose $\phi \in \mathscr S(\mathscr Q)$. Then, by $(4.2)$ and $(4.3)$, $\mathscr L\mathscr R \mathscr L(\phi)$ is well defined and
$$
\widehat{\mathscr L\mathscr R \mathscr L(\phi)}(a)= \Big(\frac{\pi}{|a|}\Big)^2 \widehat
{\phi}(a).
$$
Write $\phi_{\rho,u,v}= (\phi_{\rho})_{u,v}$. Then
$$
\mathscr L\mathscr R \mathscr L(\widetilde{\phi}_{\rho,u,v})=
\rho^{-2} \big(\widetilde{\mathscr L\mathscr R \mathscr L(\phi)}\big)_{\rho,u,v}.
$$
This equality is easy to prove by making use of the Fourier transform.

Suppose $\phi \in AW$. It is clear that
$$
W_{\phi}f=\rho^{5/2} f \ast \widetilde{\phi}_{\rho,u,v}. \eqno(5.1)
$$
Thus we can extend the wavelet transform $W_{\phi}f$ provided $(5.1)$ is valid.

Let
$$\begin{array}{rl}
L^2_{\sharp}(\mathscr Q)&=\displaystyle{\left\{ f\in L^2(\mathscr
Q): \int_{\Im \mathbb Q}\|\widehat f(a)\|^2_{HS} |a|^{-2}da<
\infty\right\},}\\[4mm]
L^2_{\natural}(\mathscr Q)&=\displaystyle{\left\{ f\in L^2(\mathscr
Q): \int_{\Im \mathbb Q}\|\widehat f(a)\|^2_{HS} |a|^{6}da<
\infty\right\}}.
\end{array}
$$
It is obvious that the Radon transform $\mathscr R$ is an
isomorphism from $L^2_{\sharp}(\mathscr Q)$ onto
$L^2_{\natural}(\mathscr Q)$. Suppose $\phi \in \mathscr S(\mathscr
Q) \bigcap AW$ and $f \in L^2_{\sharp}(\mathscr Q)$.
Then $W_{\mathscr L\mathscr R \mathscr L(\phi)}\mathscr R(f)$ is
well defined and
$$
W_{\mathscr L\mathscr R \mathscr L(\phi)}\mathscr R(f)= \rho^2 W_{\phi}f.
$$
By Theorem 5, we obtain

{\bf Theorem 10.} {\it Suppose $\phi \in \mathscr S(\mathscr Q) \bigcap AW$.
If $f \in L^2_{\sharp}(\mathscr Q)$, then
$$
 \displaystyle{f(x,t)=C_{\eta}^{-1}\int_{\mathbf
 P}W_{\mathscr L\mathscr R \mathscr L(\phi)}\mathscr R(f)(y, s, \rho, u,v)U(y, s, \rho,
 u,v)\phi(x,t)\frac {dydsd\rho du dv}{\rho^8}} \eqno(5.2)
$$
in $L^2$-sense. Furthermore, if $f\in \mathscr S(\mathscr Q) \bigcap
L^2_{_{\sharp}}(\mathscr Q)$, then $(5.2)$ holds pointwise.
Equivalently, if $f \in L^2_{\natural}(\mathscr Q)$, then
$$
 \displaystyle{\mathscr R^{-1}(f)(x,t)=C_{\eta}^{-1}\int_{\mathbf
 P}W_{\mathscr L\mathscr R \mathscr L(\phi)}f(y, s, \rho, u,v)U(y, s, \rho,
 u,v)\phi(x,t)\frac {dydsd\rho du dv}{\rho^8}} \eqno(5.3)
$$
in $L^2$-sense. Furthermore, if $f\in \mathscr S(\mathscr Q)$, then
$(5.3)$ holds pointwise.}
\vskip 0.3 true cm

If the radial wavelet $\phi_\eta$ is defined as in Section 3, then $\mathscr L\mathscr R \mathscr L(\phi)$ is also a radial function.
By Theorem 6, we have

{\bf Theorem 11.} {\it Suppose $\eta \in \mathscr S(\mathbb R_+)$ satisfies $(3.4)$ and $\phi_\eta$ is defined as in Section 3.
If $f \in L^2_{\sharp}(\mathscr Q)$, then
$$
\displaystyle{f(x,t)=C_{\eta}^{-1}\int_{\mathscr
Q\times \mathbb R_+} W_{\mathscr L \mathscr R \mathscr
L (\phi_\eta)}\mathscr R(f)(y, s, \rho)U(y, s, \rho)\phi_\eta(x,t)\frac
{dydsd\rho}{\rho^8}} \eqno(5.4)
$$
in $L^2$-sense. Furthermore, if $f\in \mathscr S(\mathscr Q) \bigcap
L^2_{_{\sharp}}(\mathscr Q)$, then $(5.4)$ holds pointwise.
Equivalently, if $f \in L^2_{\natural}(\mathscr Q)$, then
$$
\displaystyle{\mathscr R^{-1}(f)(x,t)=C_{\eta}^{-1}\int_{\mathscr
Q\times \mathbb R_+} W_{\mathscr L \mathscr R \mathscr
L (\phi_\eta )}f(y, s, \rho)U(y, s, \rho)
 \phi_\eta(x,t)\frac {dydsd\rho}{\rho^8}} \eqno(5.5)
$$
in $L^2$-sense. Furthermore, if $f\in \mathscr S(\mathscr Q)$, then
$(5.5)$ holds pointwise.}

\vskip 1 true cm

\end{document}